\documentclass{article}
\usepackage{graphicx} 
\usepackage{latexsym,color,amsmath,amsthm,amssymb,amscd,amsfonts}
\usepackage{url}
\usepackage{scalefnt}
\usepackage[font=small,labelfont=bf]{caption}
\usepackage{float}
\usepackage{graphicx}
\usepackage{amsmath}
\usepackage{amsmath}
\usepackage{relsize}
\usepackage{mathrsfs}
\usepackage{upgreek}
\usepackage{sectsty}
\usepackage{bbm}
\usepackage{tikz}
 \usetikzlibrary{patterns}

\newtheoremstyle{nonum}{}{}{\itshape}{}{\bfseries}{.}{ }{\thmnote{#3}}

\newtheorem{thm}{Theorem}[section]
\newtheorem*{thm*}{Theorem}

\newtheorem{prop}[thm]{Proposition}
\newtheorem*{prop*}{Proposition}
\newtheorem{proposition}[thm]{Proposition}

\newtheorem{cor}[thm]{Corollary}

\newtheorem{lem}[thm]{Lemma}
\newtheorem*{lem*}{Lemma}

\newtheorem{question}[thm]{Question}

\newtheorem*{rem*}{Remark}

\newtheorem*{conj*}{Conjecture}

\newtheorem*{definition*}{Definition}

\newtheorem*{fact*}{Fact}

\newtheorem*{rems*}{Remarks}

\renewcommand{\P}{\mathbb P}
\newcommand{\R}{\mathbb R}

\newcommand{\E}{\mathbb{E}}
\newcommand{\Prob}{\mathbb{P}}

\def\Vol{{\rm Vol}}

\renewcommand{\O}{\text{O}}

\def\grad{{\nabla}}

\def\eps{{\varepsilon}}

\def\conv{{\rm conv}}

\DeclareMathOperator{\dist}{dist}
\title{Combinatorial Dvoretzky Theory}
\author{Tomer Milo}
\date{}

\begin{document}

\maketitle
\begin{abstract}
    Let $P \subset \R^n $ be an $n$-dimensional polytope such that $P \subset B_2^n$. Let $M^*(P)$ denote its mean width, and let $N$ be the number of facets (faces of dimension $n-1$) of $P$. We prove:
    \begin{enumerate}
    \item If $N \leq e^{n}$, then 
    \[ M^*(P) \leq 1 - \frac{c}{\log^2(\frac{N}{n})} \]
    \item \[ M^*(P) \leq 1 - \frac{c}{n}N^{-\frac{2}{n-1}}, \]
   \end{enumerate} 
In both bounds, $c>0$ is some absolute constant.

\end{abstract}
\section{Introduction}

An $n$-dimensional polytope $P \subset \R^n$ is the convex hull of finitely many points in $\R^n$ such that $\text{int}(P)\neq \emptyset$.
$P$ has \textsl{faces} of dimensions $k = 0,1,...,n-1$. 
A general ideology in the field of asymptotic convex geometry states: 
\[ \text{If a polytope is \textsl{close} to an Euclidean ball, then it must have \textsl{many} faces.} \]
Of course, one can measure closeness between a polytope and a ball in various different ways, and furthermore, the dimension of faces whose cardinality we measure should be specified. In this paper, we shall restrict our attention to two types of faces: we denote by $V(P)$ the set $f_0(P)$ of vertices ($0$-dimensional faces) of $P$, and by $\mathcal{F}(P)$ the set $f_{n-1}(P)$ of facets (($n-1$)-dimensional faces) of $P$ (we often just write $V$ and $\mathcal{F}$ when $P$ is clear from the context). \footnote{In fact, the number of $k$-faces of $P$ is always at least $\min\{|V|,|\mathcal{F}|\}$, see \cite{Hinman2023}}

\noindent A simple example of this phenomenon is illustrated in the following proposition (which will serve as motivation for what will come next).
\begin{prop}\label{prop: close-to-ball-classic}
Let $P \subset \R^n$ be a polytope. Suppose there are numbers $0<r<R$ such that $rB_2^n \subset P \subset RB_2^n$ (here $B_2^n$ is the Euclidean unit ball in $\R^n$). Then
\[ \min\{|V|,|\mathcal{F}|\} \geq e^{\frac{1}{2}n\left(\frac{r}{R}\right)^2}
\]
Moreover, if $\frac{r}{R} = 1-\delta$ for $0 < \delta < \frac{1}{2}$, then
\[ \min\{|V|, |\mathcal{F}|\} \geq (2\delta)^{-\frac{n-1}{2}} \]
    
\end{prop}

Roughly, the argument for the lower bound on $|V|$ is as follows. Denote by $S^{n-1}$ the unit sphere, and assume without loss of generality that $R=1$. By convexity, the set $V$ must hit every $r$-cap; that is, for any $u \in S^{n-1}$ we have $V \cap \{\theta \in S^{n-1}: \langle u, \theta \rangle \geq r\} \neq \emptyset$. This implies that
\[ S^{n-1} = \bigcup_{v \in V} \{\theta \in S^{n-1}: \langle \hat{v},\theta\rangle \geq r\}  \]
where $\hat{v} := \frac{v}{|v|}$ is the normalization of $v$. Denote by $\sigma$ the unique rotationally invariant Borel probability measure on $S^{n-1}$. In the case of $r \leq \frac{1}{2}$, the $\sigma$-measure of each such cap is at most $e^{-\frac{1}{2}nr^2}$, while in the case $r = 1 - \delta > \frac{1}{2}$, the $\sigma$-measure of each cap is at most $(2\delta)^{\frac{n-1}{2}}$. The result now follows from a union bound.
The bound on $|\mathcal{F}|$ automatically follows now by applying the same argument to $P^{\circ}$ and noting that duality switches between the numbers of vertices and facets.

Informally, we would like to understand if the assumption in Proposition \ref{prop: close-to-ball-classic} about $P$ being bounded between inscribed and circumscribed balls of a given ratio can be weakened, while preserving (at least partially) the conclusion on the largeness of $|V|$ and $|\mathcal{F}|$. This paper's aim is to propose such a weakening and explore it. To this and, we take a moment to discuss a famous theorem of Dvoretzky, brought here in a formulation due to Milman (see \cite[Sections~5.2--5.3]{AGM2015} for a detailed exposition).
Every convex body $K \subset \R^n$ satisfying $0 \in \text{int}(K)$ is determined by its so called \textsl{support function}, denoted by $h_K$ which is given by
\[ h_K(x) = \max_{y \in K} \langle x,y \rangle \]
This function is homogeneous, and so is completely determined by its values on the unit sphere $S^{n-1}$. The \textsl{mean width} of $K$ is denoted by $M^*(K)$ and is defined by
\[ M^*(K) := \int_{S^{n-1}} h_K\:  d\sigma, \]
where $\sigma$ denotes the unique rotationally invariant Borel probability measure on $S^{n-1}$. As before, let $R(K) = \max_{x \in K} |x|$ denote the outradius of $K$. The following theorem is called Dvoretzky's theorem (in the formulation of Milman). See \cite[Chapter~5]{AGM2015} for an exposition.

\begin{thm}\label{thm: Dvoretzky} Let $K \subset \R^n$ be a convex body such that $0 \in \text{int}(K)$. There are absolute constants $c,c'>0$ such that, if $0<\eps<1$ and
\[
1\le k\le c\,\frac{\eps^2}{\log(2/\eps)}
\left(\frac{M^*(K)}{R(K)}\right)^2n,
\]
then a uniformly random $k$-dimensional subspace $E\subset\R^n$ satisfies
\[
(1-\eps)B_k\subset\operatorname{Proj}_E(K)\subset(1+\eps)B_k
\]
with probability at least $1-\exp(-c'\eps^2k)$. Here $\operatorname{Proj}_E$ denotes orthogonal projection onto $E$, and $B_k=M^*(K)(B_2^n\cap E)$.
\end{thm}

The \textsl{Dvoretzky parameter} $\frac{M^*(K)}{R(K)} \in (0,1]$ serves as a certain measure of the closeness of $K$ to an Euclidean ball; the bigger it is, the higher the dimension of the 'nearly Euclidean' projections $K$ admits (note that it is equal to $1$ if and only if $K$ is an Euclidean ball centered at the origin). As our ideology states, such measures are bound to have combinatorial implications when applied to polytopes.
Throughout this paper, we will assume without loss of generality that all our polytopes $P$ satisfy $R(P)=1$ and $0 \in \text{int}(P)$ (as scaling and shifting a polytope does not change its combinatorial structure). We can now state the main question this paper deals with.

\begin{question}\label{question: main}
    Let $P \subset B_2^n$ be a polytope such that $M^*(P) = 1 - \delta$. How large must $|V|$ and $|\mathcal{F}|$ be, as a function of $\delta$?
\end{question}

It is worth pausing to emphasize that using $M^*$ to measure closeness of a body $K$ to a ball is fundamentally different from the classical way of comparing the inradius and the outradius of a body. This is illustrated in the following example: Decompose orthogonally $\R^n = \R^{n-1} \times \R$. Let $B_2^{n-1} \subset \R^{n-1}$ be the Euclidean unit ball. Let $K = B_2^{n-1} \times \{0\} \subset \R^n$ be the flat $(n-1)$-dimensional unit ball embedded in $\R^n$. A standard computation reveals that $M^*(K) = 1 - \frac{1}{2n} + o(\frac{1}{n})$, yet $r=0$ and so $\frac{R}{r} = \infty$, making the flat ball very close to a full dimensional ball in $M^*$, yet as far away from a ball as possible in the `classical' sense.

It turns out that, unlike in Proposition~\ref{prop: close-to-ball-classic} (where due to duality, these questions had the same answer), these questions now have very different answers. We shall start with the bound on $|V|$ which, as it turns out, is much simpler. Actually, it is nearly identical to the answer we got in Proposition~\ref{prop: close-to-ball-classic}, as is demonstrated in the following Proposition.

\begin{prop}\label{prop: many vertices M*} Let $P \subset \mathbb{R}^{n}$ be a polytope such that $P \subset B_2^n$. Denote $M^*(P) = r$. Then
\[ |V| \geq \frac{1}{2}e^{\frac{1}{8}nr^2}.\]
Furthermore, if $r = 1 - \delta > \frac{1}{2}$, then
\[ |V| \geq \frac{1}{2}(4\delta)^{-\frac{n-1}{2}}, \]

\end{prop}

\begin{proof}
We prove the bound in the $1-\delta$ case this time. By Markov's inequality for the function $1 - h_{P}$ on $S^{n-1}$, we have

\[ \sigma(\{\theta \in S^{n-1}: h_P \geq 1 - 2\delta) \geq \frac{1}{2} \]

Denoting this `good' set by $A$, we have that for all $\theta \in A$, there exists some $v \in V$ for which $\langle v ,\theta \rangle > 1 - 2\delta$, and in particular, $\langle \hat{v}, \theta \rangle > 1- 2\delta$. A union bound gives 
\[\frac{1}{2} \leq \sigma(\bigcup_{v \in V }\{\theta \in S^{n-1}: \langle \theta, \hat{v} \rangle > 1 - 2\delta\}) \leq  \sum_{v \in V} \sigma(\{\theta \in S^{n-1}: \langle \theta, \hat{v} \rangle > 1 - 2\delta\}) \leq |V| (4\delta)^{\frac{n-1}{2}} \]
where in the $\leq$ step we have used the standard estimate on the measure of a small cap (for a which a proof can be found in, for example, \cite[Lemma~2.5]{FM2024}). This concludes the proof of the $1-\delta$ case.

 \end{proof}

The rest of the paper will be devoted to the task of bounding $|\mathcal{F}|$ in terms of $\delta$. Before we get to it, we present an example that is meant to highlight one of the difficulties of this endeavor. Let $B_2^{n-1} \times \{0\}\subset \R^n = \R^{n-1} \times \R e_n$ be the flat ($n-1$)-dimensional disk, as before. We `polytopize' it by finding a polytope $Q \subset \R^{n-1}$ such that $(1-\frac{1}{2n})B_2^{n-1} \subset Q \subset B_2^{n-1}$ ($Q $ can be constructed by taking the convex hull of any maximal subset $V \subset S^{n-2} = \partial B_2^{n-1}$ with respect to the condition $\langle v,u \rangle < 1-\frac{1}{2n}$ for all different $v,u \in V$). To transform $Q$ into a legitimate $n$-dimensional polytope, we let
$P = \conv(Q, he_n)$, where $h>0$ is some small parameter ($h \leq \frac{1}{n^2}$ suffices). Then, $P$ satisfies $M^*(P) = 1 - \frac{1}{2n} + o(\frac{1}{n})$, while all the unit normals to facets of $P$ (except for a single one!) concentrate on a very small cap around $e_n$. Recall that the idea behind the proof of Proposition~\ref{prop: many vertices M*} was that the set of normalized vertices $\hat{V}$ forms a dense net for most of $S^{n-1}$; as this example demonstrates, this can be far from the truth in the case of the unit normals to the facets. 

To summarize: we would like to say that $|\mathcal{F}|$ is large whenever $M^*(P) > 1 - \delta$ for $\delta$ small, and the last example demonstrates that if this holds, it will be despite the fact that the unit normals to $\mathcal{F}$ can be highly concentrated on $S^{n-1}$. 

Throughout this paper, the letters $c,C,c_0$ etc. denote positive absolute constants, that are not necessarily the same in different appearances.

\textbf{Use of AI tools.} This work was developed through extensive discussions with ChatGPT. These discussions were central to the development of the main proof ideas, and the model’s contribution was indispensable. I prepared the final manuscript and used ChatGPT for assistance with revisions and corrections. I take full responsibility for the mathematical content of this paper.

\textbf{Acknowledgements.} I am grateful to Amir Shpilka, Itamar Bellaïche and Rei Henigman for many useful discussions. 

\section{The few facets case}

It turns out that Proposition~\ref{prop: many vertices M*} has no similar 'dual' result for the number of facets $|\mathcal{F}|$, and that this question is fundamentally different. To get a feel for it, consider the following example. Let $P = \frac{1}{\sqrt{n}}[-1,1]^n \subset B_2^n$ be the normalized unit cube. A standard computation (see \cite[Example~3.3]{KLZ2017}) gives $M^*(P) = \sqrt{\frac{2}{\pi}} + o(1)$. By Proposition~\ref{prop: many vertices M*}, the cube has at least exponentially many vertices (as it should). However, the number of facets is only $2n$.

We will prove two bounds, useful in different ranges of the number of facets and of the deficit $\delta=\delta_n$. This section is devoted to the proof of the following theorem, which treats the first of these cases.

\begin{thm}\label{thm: first M* bound}
    Let $P \subset \R^n$ be a polytope such that $P \subset B_2^n$. Let $N = |\mathcal{F}|$ denote its number of facets, and assume $N \leq e^n$. Then,
    \[ M^*(P) \leq 1 - \frac{c}{\log^2(\frac{eN}{n})} \]
where $c >0$ is some universal constant.
\end{thm}
We conjecture that this bound is not tight, and the correct bound on $1 - M^*(P)$ should be of order $\frac{1}{\log(\frac{eN}{n})}$; see Section~\ref{sec:gaussian-tightness} for a detailed discussion on the question of tightness.

The proof of Theorem~\ref{thm: first M* bound} will require several lemmas. We start building the setup for the proof now. For a polytope $P \subset \R^n$, let $p: S^{n-1} \to P$ denote the nearest point projection:
\[ p(\theta) = \operatorname*{argmin}_{x \in P}|x-\theta| \]
Let $F_{\theta}$ be the unique face of $P$ such that $p(\theta) \in \text{relint}(F_{\theta})$. 
The first main lemma states that generally, $\dim(F_{\theta})$ is quite large.
\begin{lem}\label{lemma: high dim of p(theta)}
Let $P \subset B_2^n$ be a polytope. In the notation above, we have
\[  \int_{S^{n-1}} \dim(F_{\theta}) d\sigma \geq  (n-1)\int_{S^{n-1}} \langle p(\theta), \theta \rangle d\sigma \]
    
\end{lem}

\begin{proof}[Proof of Lemma~\ref{lemma: high dim of p(theta)}] We use the integration by parts on $S^{n-1}$: 
\begin{equation}\label{integration by parts} \int_{S^{n-1}} \langle \grad_Sf, \grad_S g \rangle d\sigma = - \int_{S^{n-1}} f\Delta_S \: g \: d\sigma, 
\end{equation}
where $\grad_S$ denotes the spherical gradient, and $\Delta_S = div_S \grad_S$. Noting that the nearest point projection to a polytope is piecewise affine (see \cite[Theorem~3.1]{FP2025}), we write $p = (p_1,\ldots,p_n)$ and apply this formula to $f = p_i$ and $g = \theta_i$. We have $\grad_S(\theta_i) = e_i - \theta_i\theta$ and so $\Delta_S \theta_i = -(n-1)\theta_i$. So. equation \ref{integration by parts} becomes
\begin{equation}\label{integration by parts 2} \int_{S^{n-1}} \langle \grad_S p_i, \grad_S \theta_i \rangle d\sigma = (n-1) \int_{S^{n-1}} p_i \theta_i d\sigma. \end{equation}
We apply equality \ref{integration by parts 2} for $i=1,\ldots,n$ and sum both sides. The right hand side becomes 
\begin{equation}\label{sum of RHS} (n-1)\int_{S^{n-1}} \langle p(\theta),\theta \rangle d\sigma. \end{equation}
We estimate the left hand side. 
\begin{equation}\label{messy}
\begin{aligned}
\int_{S^{n-1}}\langle\grad_Sp_i,\grad_S\theta_i\rangle\,d\sigma
&=\int_{S^{n-1}}\langle\grad p_i,e_i-\theta_i\theta\rangle\,d\sigma\\
&=\int_{S^{n-1}}\frac{\partial p_i}{\partial x_i}\,d\sigma
-\int_{S^{n-1}}\theta_i\bigl(Dp(\theta)\theta\bigr)_i\,d\sigma,
\end{aligned}
\end{equation}
where $Dp(\theta)$ is the Jacobian of $p$ at $\theta$, defined almost everywhere on the sphere.

Summing (\ref{messy}) over $i$ and using (\ref{integration by parts 2}) and (\ref{sum of RHS}) gives

\[
(n-1)\int_{S^{n-1}}\langle p(\theta),\theta\rangle\,d\sigma
=\int_{S^{n-1}}\operatorname{tr}Dp(\theta)\,d\sigma
-\int_{S^{n-1}}\langle\theta,Dp(\theta)\theta\rangle\,d\sigma.
\]
For almost every $\theta$, $p$ agrees locally with the affine orthogonal projection onto $\operatorname{aff}(F_\theta)$. Thus $Dp(\theta)$ is the orthogonal projection onto $H_\theta=\operatorname{span}(F_\theta-F_\theta)$, and
\[
\operatorname{tr}Dp(\theta)=\dim(F_\theta),\qquad
\langle\theta,Dp(\theta)\theta\rangle
=|\operatorname{Proj}_{H_\theta}\theta|^2\ge0.
\] Comparing (\ref{sum of RHS}) with (\ref{messy}), with the above discussion, gives
\[ \int_{S^{n-1}} \dim(F_{\theta})d\sigma \geq (n-1) \int_{S^{n-1}} \langle p(\theta), \theta \rangle d\sigma. \]
\end{proof}
\noindent The following lemma demonstrates that $\int_{S^{n-1}} \langle p(\theta), \theta \rangle d\sigma$ is large whenever $M^*(P)$ is.

\begin{lem}\label{lem: p(theta) close to theta often} Let $P \subset B_2^n$ be a polytope and let $p$ be the nearest point projection. Suppose \\ $M^*(P) > 1 - \delta$. Then,
 \[ \int_{S^{n-1}} \langle p(\theta), \theta \rangle d\sigma \geq  1 - \sqrt{2\delta}. \]
\end{lem}
\begin{proof} Let $\theta \in S^{n-1}$, and fix some $v \in V$.
    \[ |\theta - p(\theta)|^2 \leq |\theta - v|^2 = 1 - 2\langle \theta,v \rangle + |v|^2 \leq  2 - 2 \langle\theta, v \rangle. \]
Since this holds for all $v \in V$, and $h_{P}(\theta) = \displaystyle \max_{v\in V} \langle v, \theta \rangle$, we have
\[ |\theta - p(\theta)|^2 \leq  2(1-h_{P}(\theta)). \]
In particular, if $M^*(P) \geq 1 -\delta$, then
\[ \int_{S^{n-1}}|\theta - p(\theta)|^2 d\sigma \leq 2 \delta. \]
Hence $\langle\theta,p(\theta)\rangle\ge1-|p(\theta)-\theta|$, and Cauchy--Schwarz gives
\[ \int_{S^{n-1}} \langle \theta,p(\theta)\rangle d\sigma \geq 1 - \int_{S^{n-1}}|p(\theta)-\theta|d\sigma \geq 1-\sqrt{2\delta}.  \]
\end{proof}
\noindent Combining Lemma~\ref{lem: p(theta) close to theta often} and Lemma~\ref{lemma: high dim of p(theta)} gives one of the main ingredients in the proof of Theorem~\ref{thm: first M* bound}.

\begin{cor}\label{cor: main cor, integral of <p,theta> is large} Let $P \subset B_2^n$ be a polytope with $M^*(P) \geq 1 - \delta$. Let $ F_{\theta}$ be as before. Then
    \[ \int_{S^{n-1}} \dim(F_{\theta}) d\sigma > (n-1)(1 - 2\sqrt{\delta})  \]
\end{cor}

\noindent The next lemma states that the difference $\theta - p(\theta)$ is always in the normal cone of $p(\theta)$.

\begin{lem}\label{lem: normal cone}
    Let $P \subset \R^n$ be a polytope. Let $\theta \in S^{n-1}$ and let $p(\theta)$ be the nearest point to $\theta$ in $P$. Let $F_{\theta}$ be the unique face of $P$ such that $p(\theta) \in \text{relint}(F_{\theta})$. Then $q(\theta):=\theta - p(\theta) \in \text{NC}_P(p(\theta))$, that is, 
    \[ \langle \theta - p(\theta), y - p(\theta) \rangle \leq 0, \: \forall y \in P \]
\end{lem}
\begin{proof}
    Let $\phi(s) = |\theta - (p(\theta) +s(y-p(\theta)))|^2$ for $s \in [0,1]$.
    \[ \phi'(s) = -2 \langle \theta - p(\theta) - s(y - p(\theta)), y- p(\theta) \rangle \]
 By assumption, $\phi$ attains its minimum at $s=0$. Since $\phi$ is convex, $\phi_{+}^{'}(0) \geq 0$. But
 \[ 0 \leq \phi_{+}'(0) = -2 \langle \theta - p(\theta), y - p(\theta) \rangle\]
\end{proof}
We will now give an informal outline of how one utilizes this information to prove Theorem~\ref{thm: first M* bound}. Let $u_1,...,u_N$ be the unit normals to the facets of $P$. The normal cone $\text{NC}_P (p(\theta))$ is spanned by a subset $u_{i_1},\ldots,u_{i_{k(\theta)}}$, where $k(\theta) = n - \dim(F_{\theta})$. Corollary \ref{cor: main cor, integral of <p,theta> is large} tells us that $k(\theta) = \dim(E_{\theta}) \le C\sqrt{\delta}\: n$ for most $\theta$'s (after applying Markov's inequality, where $C$ is some large constant which has to be carefully chosen),
while Lemma~\ref{lem: normal cone} tells us that $q(\theta) \in E_{\theta}$. Thus, for most $\theta$'s, $q(\theta)$ must fall in at most $\binom{N}{\lceil c\sqrt{\delta}n \rceil}$ possible spaces, each of dimension at most $c\sqrt{\delta}(n-1) $. \\
\indent The final major step of the proof would be to conclude that since $q(\theta) \in E_{\theta}$, it must be that $\operatorname{dist}(\theta, E_{\theta})$ is small. For most $\theta$'s, as explained above, there are only a handful of such $E_{\theta}$'s, and they all have relatively small dimension. The proof of Theorem~\ref{thm: first M* bound} will follow by contradiction; if $\delta$ is too small relative to $N$, there are too few subspaces of too small dimension to capture most $\theta$'s in their small neighborhoods.\\
\indent The trouble with this plan is that generally, $q(\theta) \in E_{\theta}$ does not give sufficient information on $\text{dist}(\theta, E_{\theta})$, even when one has good control on $|q(\theta)|$. It turns out that in order to get that implication, one extra trick is needed, which is to run the whole argument on a slightly scaled version of $P$; One needs to work all along with $(1-\eps)P$ instead of $P$, where $\eps$ is a parameter (which will eventually be taken as a sufficiently large multiple of $\sqrt{\delta}$). The following lemma is a modification of Lemma~\ref{lem: p(theta) close to theta often} and Lemma~\ref{lemma: high dim of p(theta)}, which gives an analogue of Corollary~\ref{cor: main cor, integral of <p,theta> is large} for the scaled polytope $(1-\eps)P$.

\begin{lem}\label{lem: info on scaled polytope} Let $P \subset B_2^n$
be a polytope such that $M^*(P) \geq 1 - \delta$. Let $\eps \in (0,\frac{1}{2})$ and let $\tilde{P} = (1-\eps)P$. Let $\tilde{p}:S^{n-1} \to \tilde{P}$ be the nearest point projection, and let $\tilde{q}(\theta) = \theta - \tilde{p}(\theta)$. Let $F_{\theta}$ be the unique face for which $\tilde{p}(\theta) \in \text{relint}(F_{\theta})$. We have
\[ \int_{S^{n-1}} |\tilde{q}(\theta)-\eps \theta|^2 d\sigma \leq 2 \delta, \]
and furthermore,
\[ \int_{S^{n-1}} \dim(F_{\theta}) d\sigma \geq (1 - (\eps + \frac{\delta}{\eps}))(n-1).   \]
    
\end{lem}

\begin{proof} For any vertex $v \in V$ of $P$, we have
\[ |\tilde{q}(\theta)|^2 \leq |\theta - (1-\eps)v|^2 = \langle \theta - (1-\eps)v, \theta - (1-\eps)v \rangle \leq 1 + (1-\eps)^2 -2(1-\eps)\langle v, \theta \rangle . \]
So
\begin{equation}\label{equation q integral small} |\tilde q(\theta)|^2\le\eps^2+2(1-h_P(\theta)). \end{equation}
Moreover,
\begin{equation}\label{equation on d} |\tilde{q}(\theta)|^2 = |\tilde{q}(\theta) - \eps \theta +\eps \theta|^2 = |\tilde{q}(\theta)-\eps \theta|^2 + \eps^2 + 2\eps\langle \theta, \tilde{q}(\theta)-\eps\theta \rangle.  \end{equation}
Let $d(\theta):=\langle\theta,\tilde q(\theta)-\eps\theta\rangle$. Since $|\tilde p(\theta)|\le1-\eps$, we have $d(\theta)=1-\eps-\langle\theta,\tilde p(\theta)\rangle\ge0$. From \eqref{equation on d} and \eqref{equation q integral small} we get
\[ \int_{S^{n-1}} |\tilde{q}(\theta) - \eps \theta|^2 d\sigma \leq 2 \delta,\]
and this proves the first part of the lemma. Equations~\eqref{equation on d} and~\eqref{equation q integral small} also give
\begin{equation}\label{equation integral of d small}
    2\eps \int_{S^{n-1}} d(\theta) d\sigma \leq 2\delta.
\end{equation}
But
\[ \langle \tilde{p}(\theta), \theta \rangle = \langle \theta - \tilde{q}(\theta),\theta \rangle = \langle (1-\eps)\theta +\eps \theta - \tilde{q}(\theta), \theta \rangle = 1 - \eps - d(\theta),  \]
and so, by (\ref{equation integral of d small}),
\begin{equation}\label{equation on p for integration lemma} \int_{S^{n-1}} \langle \tilde{p}(\theta), \theta \rangle d\sigma \geq 1 - (\eps + \frac{\delta}{\eps}).    \end{equation}
Applying Lemma~\ref{lemma: high dim of p(theta)} to the polytope $\tilde{P}$ and using (\ref{equation on p for integration lemma}), we finally conclude
\[
\int_{S^{n-1}}\dim(F_{\theta}) d\sigma \geq (n-1)(1-(\eps + \frac{\delta}{\eps})).
\]
\end{proof}
We are now ready to prove Theorem~\ref{thm: first M* bound}.
\begin{proof}[Proof of Theorem~\ref{thm: first M* bound}] Let $N := |\mathcal{F}|$ and let $L := \log(\frac{eN}{n})$.
Fix an absolute constant \(0<\rho<1/2\) such that
\[
2e^2\sqrt{\rho}\le \frac12,
\]
and choose
\[
\eps=\frac{1}{32L},
\qquad
m=\left\lceil\frac{n}{4L}\right\rceil.
\]
Suppose, towards a contradiction, that
\[
\delta\le \frac{\rho^2\eps^2}{8}.
\]
Apply Lemma~\ref{lem: info on scaled polytope} with this $\eps$, and retain its notation $\tilde P,\tilde p,\tilde q,F_\theta$. Put $k(\theta):=n-\dim(F_\theta)$. Since $\tilde P\subset(1-\eps)B_2^n$, we have $k(\theta)\ge1$, and the lemma gives $\mathbb E[k(\theta)-1]\le(n-1)(\eps+\delta/\eps)$. Since $k(\theta)$ is integer-valued, the event
\(k(\theta)>m\) implies \(k(\theta)-1\ge m\).
Consequently, the second part of Lemma~\ref{lem: info on scaled polytope} and Markov's inequality give
\[
\mathbb P(\{k(\theta)>m\})
\le \frac{n}{m}
   \left(\eps+\frac{\delta}{\eps}\right)
\le 4L\eps\left(1+\frac{\rho^2}{8}\right)
=\frac18\left(1+\frac{\rho^2}{8}\right)
<\frac14.
\]
Combining the first part of Lemma~\ref{lem: info on scaled polytope} with Markov's inequality again gives
\[
\mathbb P(\left\{
|\tilde q(\theta)-\eps\theta|>\rho\eps
\right\})
\le \frac{2\delta}{\rho^2\eps^2}
\le \frac14.
\]
Thus, on a set of spherical measure greater than \(1/2\),
both
\begin{equation}\label{equation good events}
k(\theta)\le m
\qquad\text{and}\qquad
\left|\theta-\frac{\tilde{q}(\theta)}{\eps}\right|
\le \rho
\end{equation}
hold. Set $E_\theta=\operatorname{span}\mathrm{NC}_{\tilde P}(\tilde p(\theta))$. By Lemma~\ref{lem: normal cone}, $\tilde q(\theta)/\eps\in E_\theta$, and $E_\theta$ has a basis consisting of $k(\theta)$ facet normals. Thus every such $\theta$ is within distance $\rho$ of a subspace spanned by at most $m$ facet normals.

Let \(u_1,\ldots,u_N\) be the facet normals of \(P\).
Enlarging any set of fewer than \(m\) indices to one of
size \(m\), we obtain (by (\ref{equation good events})),
\begin{equation}\label{equation measure of union of neighb. of subspaces}
\frac12
<
\sigma\left(
\bigcup_{\substack{I\subset\{1,\ldots,N\}\\ |I|=m}}
\left\{
\theta\in S^{n-1}:
\operatorname{dist}
\bigl(\theta,\operatorname{span}\{u_i:i\in I\}\bigr)
\le \rho
\right\}
\right).
\end{equation}

Since \(P\) is full-dimensional, \(N\ge n+1\), so \(L>1\)
and \(m\le n/2\). Each member of this union has a measure of at most $2^n\rho^{n-m}
$ (see \cite[Section~2, equations~(2.1)--(2.2)]{Artstein2002}).
To bound the number of sets, observe that
\[
\log\binom{N}{m}
\le m\log\left(\frac{eN}{m}\right)
=mL+m\log\left(\frac{n}{m}\right)
\le \frac n4+L+\frac ne
\le 2n.
\]
The penultimate inequality uses
\(m\le n/(4L)+1\) and
\(x\log(1/x)\le 1/e\).
The last inequality follows from
\[
L\le n+1-\log n,
\]
which is where we use \(N\le e^n\).
Applying the union bound in (\ref{equation measure of union of neighb. of subspaces}) now gives
\[
\frac12
<\binom{N}{m}\,2^n\rho^{n-m}
\le e^{2n}2^n\rho^{n/2}\\
=\left(2e^2\sqrt{\rho}\right)^n
\le 2^{-n}
<\frac12,
\]
a contradiction. Hence
\[
1-M^*(P)\geq \delta>
\frac{\rho^2\eps^2}{8}
=\frac{c}{\log^2(eN/n)},
\]
where \(c>0\) is an absolute constant.

\end{proof}

\section{The many facets case}
In this section, we prove
\begin{thm}\label{thm: second M* bound}
    Let $P \subset \R^n$ be a polytope such that $P \subset B_2^n$. Let $N = |\mathcal{F}|$. Then,
    \[ M^*(P) \leq 1 - \frac{c}{n}N^{-\frac{2}{n-1}}. \]
\end{thm}
This bound is sharper than the bound presented in Theorem~\ref{thm: first M* bound} when the number of facets $N$ is at least $e^{c\sqrt{n}}$ for some constant $c>0$. 
\noindent This theorem will be obtained as an immediate corollary of the following lemma.
\begin{lem}\label{lem: main lemma spreadness of normals}
    Let $P\subset \R^n$ be a polytope such that $P \subset B_2^n$ and $M^*(P) \geq 1 - \delta$. Denote $N := |\mathcal{F}|$ and let $u_1,\ldots,u_N$ be the unit normals to the facets of $P$. Then,
    \[ \int_{S^{n-1}} \: \displaystyle \min_{i=1,\ldots,N} \: |\theta - u_i|^2 d\sigma \leq C n\delta, \]
where $C>0$ is a universal constant.
\end{lem}

In the example given in the discussion under Proposition~\ref{prop: many vertices M*}, a polytope $P \subset B_2^n$ was constructed with the following properties: $1 -M^*(P) $ is of order $\frac{1}{n}$, and the set of unit normals to the facets of $P$ was highly concentrated on $S^{n-1}$. 
Lemma~\ref{lem: main lemma spreadness of normals} demonstrates that $ \delta $ of order $\frac{1}{n}$ is the threshold for that phenomena. Indeed, if $\delta = o(\frac{1}{n})$, this lemma implies that the $u_i$'s form a dense net of most of $S^{n-1}$, from which Theorem~\ref{thm: second M* bound} will follow.

\noindent We fix some notation for the rest of this section. For any vertex $v \in V$, let \[ \Omega_v := \{\theta \in S^{n-1}: h_P(\theta) = \langle v, \theta \rangle \}. \]
\indent Equivalently, $\Omega_v$ is the intersection of $NC_P(v)$, the normal cone of $P$ at the point $v$, with $S^{n-1}$. We remark that $NC_P(v)$ is a \textsl{pointed cone}, meaning $NC_P(v)\cap(-NC_P(v))=\{0\}$, and that it is also a \textsl{polyhedral cone}, meaning that it is given by the non-negative span of finitely many vectors (we always assume this set of generators is minimal). The sets $\Omega_v$ form a partition of $S^{n-1}$ up to a null set. On each $\Omega_v$, the minimum distance to all facet normals is at most the minimum distance to the normals generating $NC_P(v)$. The following lemma deals with the structure of such cones. It roughly states that if the center of mass of a polyhedral spherical cone is close to $S^{n-1}$, then a random point in the spherical cone is expected to be close to one of its generators. 
\begin{lem}\label{lem:cone-rays}
Let $C\subset\mathbb R^n$ be a full-dimensional cone whose generated by $u_1,\ldots,u_s \in S^{n-1}$. Assume $C$ is pointed.
If $\Theta$ is uniform on $C\cap S^{n-1}$, then
\[
\mathbb E\min_i|\Theta-u_i|^2
\le cn\bigl(1-|\mathbb E\Theta|\bigr),
\]
where $c>0$ is some universal constant.
\end{lem}

\begin{proof}
Write $\mu:=|\mathbb E\Theta|$. If $\mu\le1/2$, the conclusion
is trivial. We may therefore assume
$\mu>1/2$, and put
\[
a:=\frac{\mathbb E\Theta}{\mu} \in S^{n-1}.
\]

Let $G$ be a standard Gaussian vector conditioned to belong to $C$,
and let $g$ be an unconditioned standard Gaussian vector.
Since membership in a cone depends only on direction, conditioning
preserves the independence of the radial and angular components:
$\frac{G}{|G|}$ has the distribution of $\Theta$, whereas $|G|$ has
the distribution of $|g|$. In particular,
\[
\mathbb EG=(\mathbb E|g|)\mu a,
\qquad
\mathbb E|G|^2=n.
\]
\noindent We start by showing that
\[
a\in\operatorname{int}C
\quad\text{and}\quad
\langle a,x\rangle>0
\qquad(x\in C\setminus\{0\}).
\]
Let $\nu_F$ denote the inward unit normal of a facet $F$ of $C$.
Since $\langle\nu_F,\Theta\rangle>0$ almost surely,
$\langle\nu_F,a\rangle>0$ for every $F$, proving
$a\in\operatorname{int}C$.

For the second assertion, apply the divergence theorem
coordinatewise, using \\
$\partial_j e^{-|x|^2/2}=-x_j e^{-|x|^2/2}$ and the fact that
the outward normal on $F$ is $-\nu_F$ (formally, we apply the divergence theorem on $C \cap RB_2^n$ and take $R \to \infty$). This gives
\[
\int_C x e^{-|x|^2/2}\,dx
=
\sum_F\nu_F
\int_F e^{-|x|^2/2}\,d\mathcal H^{n-1}(x).
\]
Since the integral on the left equals
$\bigl(\int_C e^{-|x|^2/2}\,dx\bigr)\mathbb EG$, we obtain
\[
a=
\sum_F
\frac{\displaystyle
\int_F e^{-|x|^2/2}\,d\mathcal H^{n-1}(x)}
{\displaystyle
\mu\,\mathbb E|g|\int_C e^{-|x|^2/2}\,dx}
\,\nu_F.
\]
Thus $a$ is a strictly positive linear combination of all the inward
facet normals. For $x\in C$, all $\langle\nu_F,x\rangle$ are
nonnegative. If they were all zero for some $x\ne0$, the entire
line $\mathbb Rx$ would satisfy the facet inequalities defining $C$,
contrary to pointedness. Hence $\langle a,x\rangle>0$ for every
nonzero $x\in C$.
\\
It follows that
\[
Q:=\{y\in a^\perp:a+y\in C\}
\]
is a bounded polytope with nonempty interior in $a^\perp$.
The hyperplane $\langle a,x\rangle=1$ meets the extreme ray
generated by $u_i$ at $u_i/\langle a,u_i\rangle$.
Consequently, the vertices of $Q$ are
\[
z_i=\frac{u_i}{\langle a,u_i\rangle}-a,
\qquad i=1,\ldots,s.
\]
Put $r=\displaystyle \min_i|z_i|$, and decompose
\[
G=Ta+Y,\qquad T=\langle G,a\rangle,\qquad Y\in a^\perp.
\]
The strict positivity proved above gives $T>0$ almost surely.
For $t>0$ and $y\in a^\perp$,
\[
ta+y\in C
\quad\Longleftrightarrow\quad
y\in tQ.
\]
Moreover, $|ta+y|^2=t^2+|y|^2$. Thus, conditionally on $T=t$,
the density of $Y$ is proportional to
$\mathbf1_{tQ}(y)e^{-|y|^2/2}$, and the vertices of $tQ$ have norm at least $tr$. By
Lemma~\ref{lem:gaussian-polytope-moment} below, applied in dimension
$n-1$,
\[
\mathbb E(|Y|^2\mid T=t)
\ge c\min\left\{\frac{t^2r^2}{n+1},1\right\}.
\]
Since $\mu>1/2$ and $\mathbb E|g|\asymp\sqrt n$,
\[
\mathbb ET=(\mathbb E|g|)\mu\ge c\sqrt n,
\qquad
\mathbb ET^2\le\mathbb E|g|^2=n.
\]
By the Paley--Zygmund inequality,
\[
\mathbb P\left\{T\ge\frac12\mathbb ET\right\}
\ge\frac{(\mathbb ET)^2}{4\mathbb ET^2}\ge c.
\]
On this event, $T\ge c\sqrt n$. Averaging the conditional estimate
over this event gives
\[
\mathbb E|Y|^2\ge c\min\{r^2,1\}.
\]
On the other hand, radial independence gives
\[
\mathbb E|Y|^2
=n\,\mathbb E\bigl(1-\langle a,\Theta\rangle^2\bigr)
\le2n\,\mathbb E\bigl(1-\langle a,\Theta\rangle\bigr)
=2n(1-\mu). \]
Consequently,
\[
\min\{r^2,1\}\le Cn(1-\mu).
\]
It remains to translate this into proximity to one of the $u_i$'s.
The definition of $z_i$ gives
\[
|z_i|^2=\frac1{\langle a,u_i\rangle^2}-1.
\]
So
\[
\min_i|a-u_i|^2
=2\left(1-\frac1{\sqrt{1+r^2}}\right)
\le2\min\{r^2,1\}
\le Cn(1-\mu).
\]
Finally, using $\mathbb E\Theta=\mu a$,
\[
\mathbb E\min_i|\Theta-u_i|^2
\le\min_i\mathbb E|\Theta-u_i|^2
=\min_i\bigl(2-2\mu\langle a,u_i\rangle\bigr)\\
=2(1-\mu)+\mu\min_i|a-u_i|^2
\le Cn(1-\mu).
\]
\end{proof}
We now turn to the proof of Lemma~\ref{lem:gaussian-polytope-moment} that was used in the proof above.

\begin{lem}\label{lem:gaussian-polytope-moment}
Let $Q\subset\mathbb R^n$ be a full-dimensional polytope whose vertices
have norm at least $r$. Let $X$ have a probability density proportional to $\mathbf{1}_{Q}(x) e^{-\frac{|x|^2}{2}}$. Then
\[ 
\mathbb{E}|X|^2 = \frac{\int_Q|x|^2e^{-|x|^2/2}\,dx}
{\int_Qe^{-|x|^2/2}\,dx}
\ge c\min\left\{\frac{r^2}{n+2},1\right\}.
\]
\end{lem}

\begin{proof}
For $S=\operatorname{conv}\{z_0,\ldots,z_n\}$, the simplex integration
formula from
\cite[Proposition~18]{BBDKV2011} gives
\[
\frac1{|S|}\int_S|x|^2\,dx
=\frac{\sum_{i=0}^n|z_i|^2+
\bigl|\sum_{i=0}^nz_i\bigr|^2}{(n+1)(n+2)}.
\]
To apply this formula, we triangulate $Q$ using its vertices, apply the previous formula to each simplex and sum the integrals. this gives
\[ |Q|^{-1}\int_Q|x|^2\,dx\ge \frac{r^2}{(n+2)}. \]
Finally, we interpolate between the uniform measure and the Gaussian measure as follows. For $0\le s\le1$, put
\[
d\mu_s(x)=\frac{\mathbf1_Q(x)e^{-s|x|^2/2}\,dx}
{\int_Qe^{-s|y|^2/2}\,dy},
\qquad A(s)=\int|x|^2\,d\mu_s(x).
\]
Each $\mu_s$ is log-concave, hence $(-1/6)$-concave. Borell's moment inequality \cite[Lemma~7.3]{AGLLOPT2012}, with $p=4$ and $r=6$, followed by Cauchy--Schwarz, gives
\[
\int|x|^4\,d\mu_s(x)\le C \left(\int|x|d\mu_s\right)^4 \leq C A(s)^2.
\]
Consequently,
\[
A'(s)=-\frac12\left(\int|x|^4\,d\mu_s(x)-A(s)^2\right)
\ge-CA(s)^2.
\]
Integrating $(\frac{1}{A})' = - (\frac{A'}{A^2}) \le C$, and using the bound on $A(0)$, yields
\[
A(1)\ge\frac{A(0)}{1+CA(0)}
\ge c\min\left\{\frac{r^2}{n+2},1\right\}.
\]
\end{proof}
\noindent We now prove Lemma~\ref{lem: main lemma spreadness of normals}
\begin{proof}[Proof of Lemma~\ref{lem: main lemma spreadness of normals}]
Let $u_1,\ldots,u_N$ be the outward unit facet normals of $P$.
For each vertex $v$, let
\[
\Omega_v=\{\theta\in S^{n-1}:h_P(\theta)=\langle v,\theta\rangle\}.
\]
These sets partition the sphere up to a null set. The normal cone at
$v$ is full-dimensional and pointed, and its extreme rays are
generated by some of $u_1,\ldots,u_s$. Let $\Theta_v$ be uniform on $\Omega_v$. Applying Lemma~\ref{lem:cone-rays} to each normal cone, and then taking the minimum over all facet normals, gives

\begin{equation*}
\int_{S^{n-1}}\min_i |\theta-u_i|^2\,d\sigma = \sum_{v} \sigma(\Omega_v) \mathbb  E[\min_i |\Theta_v -u_i|^2]
\le Cn\sum_v\ \sigma(\Omega_v)(1-| \mathbb E\Theta_v|) \end{equation*}
\begin{equation}\label{equation integral distance 1}
\le Cn\sum_v\int_{\Omega_v}
\bigl(1-h_P(\theta)\bigr)\,d\sigma
=Cn\bigl(1-M^*(P)\bigr)\le Cn\delta.
\end{equation}
For the last inequality, use $h_P(\theta)=\langle v,\theta\rangle$ on $\Omega_v$ and Cauchy--Schwarz:
\[
\frac{1}{\sigma(\Omega_v)}\int_{\Omega_v}h_P(\theta)\,d\sigma(\theta)
=\mathbb E h_P(\Theta_v)
=\langle v,\mathbb E\Theta_v\rangle
\le|\mathbb E\Theta_v|.
\]
\end{proof}
\noindent Theorem~\ref{thm: second M* bound} now follows immediately.
\begin{proof}[Proof of Theorem~\ref{thm: second M* bound}]
Put $\delta=1-M^*(P)$, and let $u_1,\ldots,u_N$ be the outward unit facet normals. The elementary cap bound
$
\sigma\{\theta:|\theta-u|\le t\}\le t^{n-1}
$
shows that, for \\ $t=(2N)^{-1/(n-1)}$, the union of the caps of radius
$t$ has measure at most $1/2$. Thus
\begin{equation}\label{equation integral distance 2}
\int_{S^{n-1}}\min_i |\theta-u_i|^2\,d\sigma
\ge\frac{t^2}{2}\ge\frac18N^{-2/(n-1)},
\end{equation}
Combining Lemma~\ref{lem: main lemma spreadness of normals} with (\ref{equation integral distance 2})  finishes the proof.
\end{proof}

\section{Tightness of the bounds}

\subsection{Tightness of Theorem~\ref{thm: first M* bound}}\label{sec:gaussian-tightness}

Assume $N\ge n$. Let $G_1,\ldots,G_N$ be independent standard Gaussian vectors in $\mathbb R^n$, and set
\[
 Q=\conv\{\pm G_1,\ldots,\pm G_N\},\qquad
 P=r(Q)Q^\circ,\qquad L=\log(N/n).
\]
Here $r(Q)$ is the inradius about the origin. Almost surely $R(P)=1$,
$P\subset B_2^n$, and $P$ has at most $2N$ facets.

\begin{proposition}\label{prop:gaussian-upper}
There are absolute constants $C,L_0,n_0>0$ such that, whenever
$n\ge n_0$ and $L\ge L_0$,
\begin{equation}\label{eq:gaussian-precise-upper}
 \mathbb{E}[1-M^*(P)]
 \le\frac{\log\log L+C}{2L}
       +C\frac{e^{-n/10}}{\sqrt n}.
\end{equation}
\end{proposition}

\begin{proof}
Write $\kappa_n=|B_2^n|$ and $\rho_Q(\theta)=\max\{t\ge0:t\theta\in Q\}$. Polar integration and H\"older's inequality give
\begin{equation}\label{eq:gaussian-jensen}
 \begin{split}
 M^*(P)
 &=r(Q)\int_{S^{n-1}}\rho_Q(\theta)^{-1}\,d\sigma(\theta)\\
 &\ge r(Q)
       \left(\int_{S^{n-1}}\rho_Q(\theta)^n\,d\sigma(\theta)\right)^{-1/n}
 =r(Q)\left(\frac{\kappa_n}{|Q|}\right)^{1/n}.
 \end{split}
\end{equation}
We therefore need a lower bound for the inradius and an upper bound
for the volume radius of $Q$. We treat these two parameters separately below.

\paragraph{The inradius.}
Let $\varphi,\Phi$ be the standard normal density and distribution
function, and put $q(t)=2(1-\Phi(t))$.
Choose $r_0>0$ by
\[
 q(r_0)=8e^{-L}\log L.
\]
The Gaussian tail bounds
\[
 \frac{2t}{1+t^2}\varphi(t)\le q(t)\le\frac2t\varphi(t)
 \qquad(t>0)
\]
imply
\begin{equation}\label{eq:gaussian-r0}
 r_0^2=2L-\log L-2\log\log L+O(1).
\end{equation}
We claim that
\begin{equation}\label{eq:gaussian-inradius}
 \Prob\{r(Q)<r_0\}\le L^{-2n}.
\end{equation}

Almost surely, each facet consists of $n$ signed sample points with
distinct indices. For a fixed such choice, the distance of their affine
hull from the origin has distribution $|g|/\sqrt n$, where $g$ is
standard normal: the sample mean is independent of the centred
configuration and has covariance $I_n/n$.
Conditional on the selected points and their affine-hull distance $t$,
a necessary condition for a facet is that all remaining Gaussian pairs
lie in the corresponding slab. Its probability is $(1-q(t))^{N-n}$.
A union bound therefore gives
\[
 \Prob\{r(Q)<r_0\}
 \le 2^n\binom Nn\sqrt{\frac{2n}{\pi}}
       \int_0^{r_0}e^{-nt^2/2}(1-q(t))^{N-n}\,dt.
\]
The integrand is increasing on $[0,r_0]$, since its logarithmic
derivative is at least
$t((N-n)q(r_0)-n)\ge0$; here we used $2\varphi(t)\ge tq(t)$.
Bounding the integral by $r_0$ times its endpoint value gives
\[
 \begin{split}
 \Prob\{r(Q)<r_0\}
 &\le C\sqrt n\,r_0
 \exp\left(n\log\frac{2eN}{n}
              -\frac{nr_0^2}{2}-(N-n)q(r_0)\right)\\
 &\le C\sqrt{nL}\,
       \exp\left(-\frac72n\log L+n\log\log L+Cn\right)
 \le L^{-2n},
 \end{split}
\]
for a sufficiently large absolute $L_0$.

\paragraph{The volume.}
We establish the uniform estimate
\begin{equation}\label{eq:gaussian-volume}
 \left(\frac{\E|Q|}{\kappa_n}\right)^{2/n}
 \le(2L-\log L+C)
       \exp\left(C\frac{e^{-n/10}}{\sqrt n}\right).
\end{equation}
Let $K_{m,n}$ denote the nonsymmetric hull of $m$ standard Gaussian
vectors in $\mathbb R^n$. The volume estimate above will follow by combining an explicit estimate of the expected number of facets of $K_{m,n+1}$ with the identity relating this number to the expected volume of $K_{m,n}$.
\cite[Theorem~1.1, equations~(2)--(3)]{KZ} gives
\begin{align}
 \frac{\E|K_{m,n}|}{\kappa_n}
 &=\frac{m!}{n!(m-n-1)!}
   \int_{\mathbb R}\Phi(t)^{m-n-1}\varphi(t)^{n+1}\,dt,
 \label{eq:kz-ordinary}\\
 \frac{\E|Q|}{\kappa_n}
 &=\frac{N!}{n!(N-n-1)!}
   \int_0^\infty(2\Phi(t)-1)^{N-n-1}
                    (2\varphi(t))^{n+1}\,dt.
 \label{eq:kz-symmetric}
\end{align}
The facet integral in \cite[equation~(2)]{BLR}, expressed using the
standard normal density, is
\[
 \E f_n(K_{m,n+1})
 =2\binom m{n+1}\sqrt{n+1}\,(2\pi)^{n/2}
   \int_{\mathbb R}\Phi(t)^{m-n-1}\varphi(t)^{n+1}\,dt.
\]
Consequently,
\[
 \frac{\E|K_{m,n}|}{\kappa_n}
 =\frac{\sqrt{n+1}}{2(2\pi)^{n/2}}\,\E f_n(K_{m,n+1}).
\]
For $\ell=\log(m/(n+1))$, \cite[Theorem~1.1]{BLR} states
\[
 \E f_n(K_{m,n+1})
 =\frac{2^{n+1}\pi^{n/2}}{\sqrt{n+1}}\,
   \ell^{n/2}
   \exp\left(-\frac{n\log\ell}{4\ell}
       +\frac{n\vartheta}{\ell}
       +O(\sqrt n\,e^{-n/10})\right),
 \qquad -34\le\vartheta\le2,
\]
provided $n+1\ge78$ and $m\ge e^e(n+1)$. The implied constant is
absolute. Hence
\[
 \begin{split}
 \left(\frac{\E|K_{m,n}|}{\kappa_n}\right)^{2/n}
 &=2\ell\exp\left(-\frac{\log\ell}{2\ell}
             +\frac{2\vartheta}{\ell}
             +O(e^{-n/10}/\sqrt n)\right)\\
 &\le(2\ell-\log\ell+C)
               \exp(Ce^{-n/10}/\sqrt n).
 \end{split}
\]
For the symmetric hull, take $m=2N-n-1$ and use
$2\Phi(t)-1\le\Phi(t)^2$ for $t\ge0$ in
\eqref{eq:kz-symmetric}. Comparison with \eqref{eq:kz-ordinary} yields
\[
 \E|Q|
 \le\left(\prod_{j=0}^n\frac{2(N-j)}{m-j}\right)\E|K_{m,n}|
 \le e^{Cn^2/N}\E|K_{m,n}|.
\]
Taking the power $2/n$, and using
$\ell=L+\log2+O(1/n+e^{-L})$ and $Le^{-L}=O(1)$, proves
\eqref{eq:gaussian-volume}.

We combine the estimates above.
Let $A=\{r(Q)\ge r_0\}$. Applying \eqref{eq:gaussian-jensen} on $A$
and then conditional Jensen gives
\[
 \E M^*(P)
 \ge r_0\Prob(A)
       \left(\E\left[\frac{|Q|}{\kappa_n}\,\middle|\,A\right]\right)^{-1/n}
 \ge
 \frac{r_0\Prob(A)^{1+1/n}}{(\E|Q|/\kappa_n)^{1/n}}.
\]
Using \eqref{eq:gaussian-r0}--\eqref{eq:gaussian-volume}, we obtain
\[
 \E[1-M^*(P)]
 \le 1-\frac{r_0}{\sqrt{2L-\log L+C}}
          +C L^{-2n}+C\frac{e^{-n/10}}{\sqrt n}
 \le\frac{\log\log L+C}{2L}
          +C\frac{e^{-n/10}}{\sqrt n},
\]
which is \eqref{eq:gaussian-precise-upper}.

\end{proof}

\textbf{Remark.}
    In the range $L\ge L_0$ and $2N\le e^n$, there is a gap between Theorem~\ref{thm: first M* bound} and Proposition~\ref{prop:gaussian-upper}. Recall that $L=\log(N/n)$ and that $P$ has at most $2N$ facets. The theorem gives $1-M^*(P)\ge c/L^2$, whereas the proposition gives
\[
\mathbb E[1-M^*(P)]\le C\frac{1+\log\log L}{L}.
\]
In this range, the exponentially small term in \eqref{eq:gaussian-precise-upper} is absorbed into $C/L$, and some realization of $P$ satisfies the same upper bound. We conjecture that the lower bound can be improved to order $1/L$ and that a construction attains this order, removing the $\log\log L$ loss. Such an improvement will give a `continuous' transition between the bounds in Theorem~\ref{thm: first M* bound} and Theorem~\ref{thm: second M* bound}, with $N = |\mathcal{F}|$ of exponential growth being the threshold. As of now, the bound in Theorem~\ref{thm: second M* bound} becomes superior when $N$ is of order $e^{c\sqrt{n}}$, and remains of order $1 - \frac{1}{n}$ for $N \leq e^{cn}$.

\subsection{Tightness of Theorem~\ref{thm: second M* bound}}\label{sec:covering-tightness}

In this section we prove that Theorem~\ref{thm: second M* bound} is tight up to a $\log$ factor.

\begin{thm}\label{thm:construction}
There are absolute constants $C_0,C_1$ with the following properties. If $N\ge \exp(C_0n)$, there exists a full-dimensional polytope
 $P\subset B_2^n$ with at most $N$ facets, $R(P)=1$, and
 \begin{equation}\label{eq:cover-upper}
  1-M^*(P)\le C_1\frac{\log(n)}{n}N^{-2/(n-1)}.
 \end{equation}
\end{thm}

\begin{proof}[Proof of Theorem~\ref{thm:construction}].
Write $d=n-1$, and let
\[
 F_n(z)=\sigma\{u\in S^{n-1}:\langle u,e_1\rangle\ge1-z\},
 \qquad 0<z<1.
\]
For $n\ge3$ and $0<z\le1/4$, we claim that
\begin{equation}\label{eq:cap-derivative}
\frac{d+1}{3z}
\le\frac{F_n'(z)}{F_n(z)}\le\frac{d}{2z}.
\end{equation}
Denote $c_n= \frac{|S^{n-2}|}{|S^{n-1}|}$. We have
\[
 F_n(z)=c_n\int_{1-z}^1(1-t^2)^{(d-2)/2}\,dt
       =c_n\int_0^z[s(2-s)]^{(d-2)/2}\,ds.
\]
In particular, $F_n'(z)=c_n[z(2-z)]^{(d-2)/2}$, so
\[
 \frac{F_n''(z)}{F_n'(z)}
 =\frac{d-2}{2}\left(\frac1z-\frac1{2-z}\right)
 =\frac{d-2}{z}\frac{1-z}{2-z}.
\]
To compare $zF_n'(z)$ with $F_n(z)$, we compare their derivatives:
\[
 (zF_n'(z))'
 =F_n'(zs)\left(1+(d-2)\frac{1-z}{2-z}\right).
\]
Since $1/3\le(1-z)/(2-z)\le1/2$ on $0<z\le1/4$, this gives
\[
 \frac{d+1}{3}F_n'(z)\le (zF_n'(z))'\le\frac d2F_n'(z).
\]
Integrating from $0$ to $z$, using $F_n(0)=0$ and
$\lim_{s\downarrow0}sF_n'(s)=0$, yields
\[
 \frac{d+1}{3}F_n(z)\le zF_n'(z)\le\frac d2F_n(z).
\]
Dividing by $zF_n(z)>0$ proves \eqref{eq:cap-derivative}. The cap formula also gives, for $0<z\le1/4$,
\begin{equation}\label{eq:cap-lower}
F_n(z)\ge\frac{2c_n}{d}z^{d/2}\ge c n^{-1/2}z^{d/2},
\end{equation}
since $c_n\asymp\sqrt n$; choosing $C_0>0$ to be a sufficiently large absolute constant ensures that, whenever $N\ge\exp(C_0n)$, the cap heights $z_0,z_1$ defined below by $F_n(z_0)=1/N$ and $F_n(z_1)=D_n/N$, with $1<D_n\le Cn\log n$, satisfy $0<z_0<z_1\le1/4$.

For unit vectors $u_1,\ldots,u_k$, $k\le N$, put $Q=\conv\{u_i\}$, and choose $z_0\in(0,1/4]$ by $F_n(z_0)=1/N$. Since $h_Q(\theta)=\max_i\langle u_i,\theta\rangle$, the union bound gives $\sigma(\{h_Q\ge1-z\})\le NF_n(z)$. Integrating $\sigma(\{1-h_Q>z\})$ over $0<z<z_0$ therefore gives
\begin{align}
 1-M^*(Q)
 &\ge\int_0^{z_0}\bigl(1-NF_n(z)\bigr)\,dz \notag\\
 &\ge z_0\left(1-\frac{3}{d+3}\right).
 \label{eq:quantization}
\end{align}
The final inequality uses 
$F_n(z)/F_n(z_0)\le(z/z_0)^{d/3}$, which is a consequence of \eqref{eq:cap-derivative}. 

B\"or\"oczky--Wintsche's spherical covering theorem \cite{BW2003}, in the form recorded
as Theorem 2.2 of Nasz\'odi \cite{Naszodi}, states that for any $0<r<\frac{\pi}{2}$, $S^d$ admits a
covering by congruent spherical caps of angular radius $r$ and density at most
\begin{equation}\label{eq:cover-density}
 d\log d+d\log\log d+5d.
\end{equation}
Choose $1<D_n\le Cn\log n$ above this bound, and choose $z_1\in(0,1/4]$ by $F_n(z_1)=D_n/N$. Apply the covering theorem with angular radius $\arccos(1-z_1)$. If the covering has $k$ centers, its density is $kF_n(z_1)\le D_n$, so $k\le N$. For the convex hull $Q$ of these centers, the covering property gives $h_Q(\theta)\ge1-z_1$ for every $\theta$, and hence $r(Q)\ge1-z_1$. Moreover,
\begin{equation}\label{eq:quantile-gap}
 z_1-z_0\le\frac{3\log D_n}{n}\,z_1,
 \qquad z_1\le C N^{-2/d}.
\end{equation}
For the first inequality, \eqref{eq:cap-derivative} gives
\[
\log D_n=\log\frac{F_n(z_1)}{F_n(z_0)}
\ge\frac n3\log\frac{z_1}{z_0}.
\]
Thus $z_1-z_0\le z_1(1-D_n^{-3/n})\le3z_1\log D_n/n$. For the second, \eqref{eq:cap-lower} gives
\[
z_1\le\left(\frac{C\sqrt n\,D_n}{N}\right)^{2/d}
\le C N^{-2/d},
\]
because $D_n\le Cn\log n$ and $d=n-1\ge2$.

Set $P=r(Q)Q^\circ$. The elementary inequalities
$h_{Q^\circ}\ge1/h_Q$ and Jensen's inequality imply
\begin{equation}\label{eq:polar-jensen}
 M^*(P)\ge\frac{r(Q)}{M^*(Q)}.
\end{equation}
Combining \eqref{eq:quantization}--\eqref{eq:polar-jensen}, and using $z_0\le z_1\le1/4$, gives
\[
\begin{aligned}
1-M^*(P)
&\le\frac{z_1-z_0+3z_0/(d+3)}{1-z_0+3z_0/(d+3)}\\
&\le C\frac{1+\log D_n}{n}z_1
\le C\frac{\log n}{n}N^{-2/d}.
\end{aligned}
\]
This proves \eqref{eq:cover-upper}.
\end{proof}

\textbf{Remark.} We conjecture that the extra $\log n$ factor in Theorem~\ref{thm:construction} can be removed, that is, that Theorem~\ref{thm: second M* bound} is asymptotically tight.

\end{document}